\documentclass[11pt]{article}

\usepackage[utf8]{inputenc}
\usepackage[T1]{fontenc}
\usepackage{amsmath,amssymb,amsfonts}
\usepackage{graphicx}
\usepackage{bm}
\usepackage{xcolor}
\usepackage{tabularx}
\usepackage{makecell}
\usepackage{dcolumn}
\usepackage[a4paper,margin=2.5cm]{geometry}

\newcommand\pder[2][]{\ensuremath{\frac{\partial#1}{\partial#2}}}
\newtheorem{obs}{Remark}[section]
\DeclareMathOperator{\erf}{erf}

\title{Lie symmetry method for a nonlinear heat-diffusion equation}
\author{Julieta Bollati \and Ernesto A. Borrego Rodriguez \and Adriana C. Briozzo}
\date{\today}

\begin{document}
\maketitle

\begin{center}
\small Departamento de Matemática, Universidad Austral, Paraguay 1950, Rosario, 2000, Santa Fe, Argentina\\
CONICET, Argentina
\end{center}

\begin{abstract}
 We investigate the nonlinear heat-diffusion  equation
$
C(u) u_t = (K(u) u_x)_x,
$
where $C(u)$ and $K(u)$ are  coefficients that depend on $u$. By applying the classical Lie symmetry method, we determine the admitted Lie point symmetries and compute the corresponding infinitesimal generators according to the functional relationship between $C(u)$ and $K(u)$. The admitted symmetries are used to reduce the partial differential equation to ordinary differential equations and to construct invariant solutions. Particular cases of physical interest are analyzed in detail, including Storm-type materials and the dependence of power-laws of $C(u)$ and $K(u)$ on $u$, whose similarity solutions are obtained. Using a Kirchhoff-type transformation, we recover the corresponding results available in the literature for the standard nonlinear diffusion equation. We also discuss the importance of keeping $C(u)$ and $K(u)$ separated in order to preserve their physical meaning.
\end{abstract}

\section{Introduction}

The analysis of partial differential equations (PDEs) plays a central role in the formulation and interpretation of mathematical models that describe physical, biological, and technological phenomena. Such equations govern a wide range of processes, including heat and mass transfer, diffusive transport, wave propagation, fluid motion, and dynamics in continuous media. Despite their fundamental importance, obtaining exact solutions for nonlinear PDEs or for PDEs with variable or solution‑dependent coefficients remains a highly nontrivial task. Even for classical models that have been extensively investigated, the analytical difficulties arising from nonlinearity and coefficient variability impose significant restrictions on conventional solution methods. In this context, symmetry ‐-based techniques have emerged as robust and systematic tools to achieve structural simplification and analytical tractability.

Lie group theory, originally developed by Sophus Lie in the late nineteenth century \cite{Li1881}, provides a rigorous and unifying framework for the study of the invariance properties of differential equations in continuous transformation groups. The central idea is to identify the set of point transformations that leave the equation invariant and, on this basis, to construct invariant solutions, reduce the order of the equation, or decrease the number of independent variables through established reduction procedures. This methodology enables the transformation of a PDE into an ordinary differential equation (ODE) or into a PDE of reduced complexity, thereby facilitating both qualitative analysis and the derivation of explicit solutions.

Partial differential equations of heat–diffusion type constitute one of the most relevant classes of models for describing transport processes in continuous media. Such equations arise in numerous physical and applied contexts, including heat conduction, pollutant dispersion, population dynamics, flow through porous media, and various mass transfer processes \cite{CaJa1959,Cr1956}. In their simplest form, linear diffusion models involve constant coefficients and admit a classical analytical treatment. However, in many realistic situations the diffusive behavior depends on the state of the system, leading to equations with variable coefficients that depend on the unknown function. This introduces essential nonlinearities and significantly complicates the study of exact solutions.
Among these models, a particularly interesting family corresponds to nonlinear diffusion equations in which the diffusion coefficient depends explicitly on the dependent variable. Such models make it possible to describe media whose heat capacity, conductivity or diffusivity varies with concentration, temperature, or other relevant quantities.  In this context, the Lie symmetry method offers a systematic and powerful tool for the analysis of nonlinear partial differential equations. In \cite{BlKu1980, BlKu1989} had been shown that the nonlinear diffusion equation $u_t-\left(\mathcal{D}(u)u_x\right)_x=0$ with variable diffusivity $\mathcal{D}(u)$ is invariant under the three, four or five parameters point Lie groups according to  arbitrary or particular diffusivities $ \mathcal{D}(u)$. A renewed group classification of this general nonlinear heat equation was carried out in \cite{Si2025}, where the determining equations were explicitly derived and the admissible forms of $\mathcal{D}(u)$ were classified, clarifying the classical results of \cite{Ov1959}. Let's remember that a six-parameter point Lie group leaves invariant this heat  equation in the case $\mathcal{D}(u)=$ constant\cite{Li1881}.

 The classical group classification of the nonlinear heat-diffusion equation was initiated in \cite{Ov1959}, where the admissible symmetry structures were identified. The analysis was later extended in \cite{Cl1993} to the nonlinear heat equation with source term, $u_t = u_{xx} + f(u)$. Further developments in the symmetry analysis of nonlinear diffusion--convection equations $u_t=(\mathcal{D}(u)u_x)_x-K'(u)u_x$ were presented in \cite{Yu1994}, where Lie group methods were used to obtain classifications and invariant reductions.

More recently, Lie symmetry techniques have been applied to generalized heat-type models. In \cite{Al2024}, the symmetry structure and invariant solutions of a generalized modified heat equation $u_t+(u+\alpha)u_{xx}+\beta u_{xxx}=0$ were investigated.

Similar symmetry-based approaches have also been developed for other nonlinear evolution equations, including the Burgers equation and its generalizations, where Lie group techniques have been used to obtain invariant reductions and exact solutions \cite{LiLiXh2009,PeFeTz2024,Su2019,Te2010,Va2015}.

R. Cherniha has contributed extensively to the application of Lie symmetry methods to obtain exact solutions of differential equations and systems of equations \cite{ChDa2022,ChKo2012,ChSePr2021}.

The purpose of this work is threefold. We perform a Lie point symmetry analysis of a class of nonlinear heat-diffusion equations $ C(u)u_t = (K(u)u_x)_x$  with coefficients depending on the unknown function $u$. In particular, we derive the associated infinitesimal generators and determine the functional forms that the constitutive functions $K(u)$ and $C(u)$ must assume in order for the equation to admit  Lie point symmetries. We then use the resulting symmetry structure to construct invariant reductions and obtain exact and self-similar solutions. Finally, we investigate relevant particular cases previously reported in the literature. The analysis illustrates how the Lie symmetry framework provides a systematic approach to derive similarity solution families in nonlinear diffusion models.

We begin in Section 2 by performing a complete Lie point symmetry analysis for the nonlinear heat-diffusion equation $C(u)u_t = (K(u)u_x)_x$. We derive the determining equations for the infinitesimal generators and classify the admitted symmetries according to the functional relationship between $C(u)$ and $K(u)$. In Section 3, we apply these symmetries to construct invariant solutions. For each admitted generator, we reduce the original partial differential equation to ordinary differential equations or simpler forms, obtaining families of similarity solutions. In Section 4, we apply the theoretical results to three particular cases found in the literature. For these models, we derive explicit or implicit invariant solutions,  extending previously known results. In Section 5, the results obtained for the present model are compared with those derived for the standard nonlinear diffusion equation studied in \cite{BlKu1989, Ov1959}. The correspondence between the infinitesimal generators of both equations is established, and the relationships between their respective coefficients are analyzed. Particular attention is devoted to the importance of preserving the explicit separation of the coefficients $K$ and $C$, highlighting the physical information that may be lost when both quantities are combined into a single effective diffusion coefficient.
Finally, Section 6 presents the main conclusions of this work and discusses potential avenues for future research, such as the inclusion of source terms or the application to free boundary problems.

\section{Lie point symmetries}
Let us consider the following nonlinear partial differential equation, commonly referred to as the heat-diffusion equation, in which we assume that the coefficients $C=C(u)$ and $K=K(u)$ are non zero and non simultaneously constant
\begin{equation} \label {ec}
    C(u)u_t = (K(u)u_x)_x,
\end{equation}
which can be written as
\begin{equation}
    C(u)u_t = K'(u)u_x^2 + K(u)u_{xx}.\label{eq}
\end{equation}

Our goal is to determine the one-parameter Lie group of  transformations depending on the parameter $\epsilon$ given by
\begin{equation}\label{groupLie}
\left\{
\begin{array}{l}
    x^{*}=R(x,t,u,\epsilon)\\
    t^{*}=T(x,t,u,\epsilon)\\
    u^{*}=U(x,t,u,\epsilon),
    \end{array}
    \right.
\end{equation}
that map the space $(x, t, u)$ into itself, under which equation \eqref{ec} remains invariant.
The infinitesimal transformation of the one-parameter Lie group of transformations \eqref{groupLie} is given by
\begin{equation}\label{group}
\left\{
\begin{array}{l}
    x^{*}=x+\epsilon \xi_1(x,t,u)+o(\epsilon^{2})\\
    t^{*}=t+\epsilon \xi_2(x,t,u)+o(\epsilon^{2})\\
    u^{*}=u+\epsilon \eta(x,t,u)+o(\epsilon^{2}),
    \end{array}
    \right.
\end{equation}
where 
\begin{equation}
\left\{
\begin{array}{l}
    \xi_1(x,t,u)=\pder[R]{\epsilon}\big|_{\epsilon=0}\\[6pt]
    \xi_2(x,t,u)=\pder[T]{\epsilon}\big|_{\epsilon=0}\\[6pt]
    \eta(x,t,u)=\pder[U]{\epsilon}\big|_{\epsilon=0}
    \end{array}
    \right.
\end{equation}
are called the infinitesimals of \eqref{groupLie}.

If we denote
$$
F(u,u_x,u_t,u_{xx}) = C(u)u_t - K'(u)u_x^2 - K(u)u_{xx},
$$
then the invariance of equation \eqref{ec} under the transformation group \eqref{group} is equivalent to the invariance condition \cite{BlKu1989,Hy2000, Ol1986, Ov1982}
\begin{equation}\label{X2F=0}
X^{(2)}(F) = 0 \text{ when } F=0,
\end{equation}
where $X^{(2)}$ is the second prolongation of the infinitesimal generator \begin{equation}
X=\xi_1\frac{\partial}{\partial x} 
+ \xi_2\frac{\partial}{\partial t} 
+ \eta\frac{\partial}{\partial u}
\end{equation}  which is given by 
\begin{align}
X^{(2)} = X+ \eta_1^{(1)}\tfrac{\partial}{\partial u_x} 
+ \eta_2^{(1)}\tfrac{\partial}{\partial u_t}
+ \eta_{11}^{(2)}\tfrac{\partial}{\partial u_{xx}}
+ \eta_{12}^{(2)}\tfrac{\partial}{\partial u_{xt}}
+ \eta_{21}^{(2)}\tfrac{\partial}{\partial u_{tx}}
+ \eta_{22}^{(2)}\tfrac{\partial}{\partial u_{tt}}.
\end{align}
Since $\pder[F]{x}=\pder[F]{t}=\frac{\partial F}{\partial u_{xt}}=\frac{\partial F}{\partial u_{tx}} = \frac{\partial F}{\partial u_{tt}}=0$, we have
\begin{equation}\label{X2F}
\begin{aligned}
X^{(2)}(F) =  
 \eta\frac{\partial F}{\partial u}
+ \eta_1^{(1)}\frac{\partial F}{\partial u_x} 
 + \eta_2^{(1)}\frac{\partial F}{\partial u_t}
+ \eta_{11}^{(2)}\frac{\partial F}{\partial u_{xx}},
\end{aligned}
\end{equation}
where
\begin{equation*}
\eta_{1}^{(1)} = \frac{\partial \eta}{\partial x} 
+ \left[ \frac{\partial \eta}{\partial u} 
- \frac{\partial \xi_{1}}{\partial x} \right] u_{x} 
- \frac{\partial \xi_{2}}{\partial x} u_{t} - \frac{\partial \xi_1}{\partial u}u_x^2 - \frac{\partial \xi_2}{\partial u}u_xu_t,
\end{equation*}
\begin{equation*}
\eta_{2}^{(1)} = \frac{\partial \eta}{\partial t} 
+ \left[ \frac{\partial \eta}{\partial u} 
- \frac{\partial \xi_{2}}{\partial t} \right] u_{t} 
- \frac{\partial \xi_{1}}{\partial t} u_{x} - \frac{\partial \xi_2}{\partial u}u_t^2 - \frac{\partial \xi_1}{\partial u}u_xu_t,
\end{equation*}
\begin{equation*}
\begin{aligned}
\eta_{11}^{(2)} &= \frac{\partial^{2} \eta}{\partial x^{2}} 
+ \left[ 2 \frac{\partial^{2} \eta}{\partial x \partial u} 
- \frac{\partial^{2} \xi_{1}}{\partial x^{2}} \right] u_{x}
- \frac{\partial^{2} \xi_{2}}{\partial x^{2}} u_{t} + \left[ \frac{\partial \eta}{\partial u} - 2 \frac{\partial \xi_1}{\partial x} \right]u_{xx} \\[6pt]
&\quad  -2\frac{\partial \xi_2}{\partial x}u_{xt} + \left[ \frac{\partial^{2} \eta}{\partial u^{2}} - 2 \frac{\partial^{2} \xi_{1}}{\partial x \partial u} \right] u_{x}^{2}
- 2 \frac{\partial^{2} \xi_{2}}{\partial x \partial u} u_{x} u_{t} - \frac{\partial^2 \xi_1}{\partial u^2}u_x^3 \\[6pt]
&\quad  - \frac{\partial^2 \xi_2}{\partial u^2}u_x^2u_t - 3\frac{\partial \xi_{1}}{\partial u} u_{x} u_{xx}
- \frac{\partial \xi_{2}}{\partial u} u_{t}u_{xx} - 2\frac{\partial\xi_2}{\partial u}u_xu_{xt}.
\end{aligned}
\end{equation*}
are the extended infinitesimals.
The invariance condition \eqref{X2F=0} applied to equation \eqref{eq} leads to 
\begin{equation}\label{invcond}
\eta\left(C'(u) u_t-K''(u) u_x^2-K'(u) u_{xx}\right)-2K'(u) u_x\,\eta_1^{(1)} 
+\eta_2^{(1)} C(u)-\eta_{11}^{(2)} K(u)=0 .
\end{equation}
By using \eqref{eq} to eliminate $u_t$, and collecting the coefficients of $u_{xt}$, $u_xu_{xt}$, $u_xu_{xx}$ in \eqref{invcond}, we obtain
\begin{equation}
    \frac{\partial \xi_1}{\partial u} = \frac{\partial \xi_2}{\partial x} = \frac{\partial \xi_2}{\partial u} = 0,
\end{equation}
which implies that 
\begin{align}
\xi_1 = \xi_1(x,t) ,\qquad\;\xi_2 = \xi_2(t).
\end{align}

Terms that do not involve any derivative of $u$ lead to
\begin{equation}
    C(u)\frac{\partial \eta}{\partial t} - K(u)\frac{\partial^2 \eta}{\partial x^2} = 0, \label{ti}
\end{equation}
and from the coefficients $u_x$, $u_x^2$ $u_{xx}$ in \eqref{invcond} we have
\begin{equation}
    2K'(u)\frac{\partial \eta}{\partial x} 
    + C(u)\frac{\partial \xi_1}{\partial t} 
    + K(u)\left[ 2\frac{\partial^2 \eta}{\partial x \partial u} - \frac{\partial^2 \xi_1}{\partial x^2}\right] = 0, \label{ux}
\end{equation}
\begin{equation}
    \eta\left[ \frac{C'(u)K'(u)}{C(u)} - K''(u) \right] 
    + K'(u)\left[ 2\frac{\partial \xi_1}{\partial x} - \frac{\partial \eta}{\partial u} - \xi_2'\right] 
    - K(u)\frac{\partial^2 \eta}{\partial u^2} = 0, \label{ux2}
\end{equation}
\begin{equation}
    \eta\left[ \frac{C'(u)K(u)}{C(u)} - K'(u) \right] 
    + K(u)\left[ 2\frac{\partial \xi_1}{\partial x} - \xi_2' \right] = 0. \label{uxx}
\end{equation}

Solving \eqref{uxx} for $\eta$, we obtain
\begin{equation}\label{eta}
    \eta = A(u)\left( \xi_2'- 2\frac{\partial \xi_1}{\partial x} \right),
\end{equation}
where
\begin{equation}\label{A}
A(u) = \left( \tfrac{C'(u)}{C(u)} - \tfrac{K'(u)}{K(u)}\right)^{-1},
\end{equation}
and substituting \eqref{eta} into \eqref{ti}, we obtain
\begin{equation}\label{casos}
    \frac{C(u)}{K(u)}\left[ 2\frac{\partial^2 \xi_1}{\partial x \partial t} - \xi_2'' \right] 
    = 2\frac{\partial^3 \xi_1}{\partial x^3}.
\end{equation}

The factor $\tfrac{C(u)}{K(u)}$ depends only on the variable $u$, while the remaining terms depend only on the independent variables $x$ and $t$. This observation implies that the equation \eqref{casos} can only be satisfied under different structural conditions. Therefore, the analysis naturally separates into two distinct cases, thus, two cases can be distinguished: the ratio
$
\frac{C(u)}{K(u)}
$
is either non-constant, or \(\frac{C(u)}{K(u)} = \alpha\) constant.

Next we consider each case to obtain the corresponding one-parameter Lie group of transformation \eqref{groupLie}.

\subsection*{Case 1:  $\tfrac{C(u)}{K(u)}$ non constant }
As in this case $\frac{C(u)}{K(u)}$ depend on $u$, it follows from \eqref{casos} that
\begin{equation}\label{xi13=0}
\pder[^3\xi_1]{x^3} = 0, \qquad 
 2 \pder[^2\xi_1]{x \partial t}-\xi_2''  = 0, 
\end{equation}
then we find 
\begin{equation}\label{xi1}
\xi_1(x,t) = \rho x^2 + \mu(t) x + \gamma(t), \qquad \xi_2''(t) =2\mu'(t)
\end{equation}
where $\rho$ is a constant, and $\mu$ and $\gamma$ are arbitrary functions of $t$. As a consequence,
\begin{equation}
    \xi_2'(t)=2\mu(t)-2\beta,
\end{equation}
where $\beta$ is a  constant of integration and hence
\begin{equation}\label{xi1-f}
    \xi_1(x,t)=\rho x^2+\left( \frac{\xi_2'(t)}{2}+\beta\right) x+\gamma(t).
\end{equation}
From \eqref{eta}, it immediately follows that
\begin{equation}\label{eta-f}
    \eta(x,t,u)=-A(u)\left(4\rho x+2\beta \right).
\end{equation}
Taking into account  \eqref{ux}, we obtain that 
\begin{equation}
    -8\rho \frac{K'(u)}{C(u)}A(u)+\frac{K(u)}{C(u)} \left( -8\rho A'(u)-2\rho\right)+ \frac{\xi_2''(t)}{2}x+\gamma'(t)=0,
\end{equation}
and then separating the terms that depend on $u$ from those that do not, we get
\begin{equation}\label{rhoKA}
    2 \rho\left[4 \frac{K'(u)}{C(u)}A(u)+\frac{K(u)}{C(u)} \left( 4 A'(u)+1\right)\right]=0
\end{equation}
and
\begin{equation}
    \frac{\xi_2''(t)}{2}x+\gamma'(t)=0. \label{xi2gamma}
\end{equation}
From \eqref{rhoKA}, we deduce  that 
\begin{equation}
    \rho=0\qquad   \text{or} \qquad 4 \frac{K'(u)}{C(u)}A(u)+\frac{K(u)}{C(u)} \left( 4 A'(u)+1\right)=0,
\end{equation}
and from \eqref{xi2gamma}, it follows that $\gamma(t)\equiv\gamma$ is  constant and $\xi_2$ takes the form
\begin{equation} \label{xi2}
    \xi_2(t)=\delta t+\sigma,
\end{equation}
where $\delta,\sigma,\gamma$ are arbitrary constants.

Consequently, from \eqref{xi1-f}, \eqref{eta-f} and \eqref{xi2},  the infinitesimals are given by
\begin{equation}
\left\{
\begin{array}{l}
\xi_1(x) = \rho x^2 + \left( \frac{\delta}{2} + \beta \right) x + \gamma, \\
\xi_2(t) = \delta t + \sigma, \\
\eta(x,u) = -A(u) \left( 4 \rho x + 2 \beta \right),
\end{array}\right.
\end{equation}
with  five parameters $\delta, \beta,\rho, \gamma, \sigma$.

It remains to verify that equation \eqref{ux2}  holds, which is equivalent to 
\begin{equation}\label{ux2-1}
(4\rho x+2\beta)\Big[ 
 -A(u)\Big(\tfrac{C'(u)}{C(u)}K'(u)-K''(u)\Big)
 +K'(u) 
 +K'(u)A'(u)+K(u)A''(u)
\Big]=0 .
\end{equation}
 Notice that the parameters $\delta$, $\sigma$ and $\gamma$ exist for arbitrary $K(u)$ and $C(u)$ but from \eqref{ux2-1}, the existence of parameters $\beta$ and $\rho$ depends on the form of $K(u)$ and $C(u)$. Three cases arise:

\medskip

\noindent \textbf{Case 1.1:} $\rho = 0$, $\beta = 0$.  \\
In this case, equations \eqref{rhoKA} and \eqref{ux2-1} are automatically satisfied for any functions $K(u)$ and $C(u)$ such that $C(u)/K(u)$ is non-constant.  
Consequently, the PDE \eqref{eq} admits a three-parameter Lie group with infinitesimal generators
\begin{equation}\label{infgen}
X_1 = \frac{x}{2} \pder[]{x} + t \pder[]{t}, \qquad
X_2 = \pder[]{x}, \qquad
X_3 = \pder[]{t},
\end{equation}
corresponding to the parameters $\delta$, $\gamma$ and $\sigma$ respectively.

\medskip
\noindent \textbf{Case 1.2:} $\rho = 0$, $\beta\neq 0$.\\
In this case, equation \eqref{rhoKA} is automatically satisfied, whereas equation \eqref{ux2-1} leads to the following condition
\begin{equation}
A(u) \left( \tfrac{K'(u)}{K(u)}\right)'+\left( \tfrac{K'(u)}{K(u)}\right) A'(u)+A''(u)=0,
\end{equation}
which can be equivalently rewritten as
\begin{equation}
    \left(A(u) \tfrac{K'(u)}{K(u)}\right)'=-A''(u).
\end{equation}
By integrating in $u$, one finds
\begin{equation}
    \left(A(u)K(u)\right)'=BK(u),
\end{equation}
where $B$ is an arbitrary constant.

As a consequence we obtain the following expression for $A$
\begin{equation}\label{Aint}
    A(u)=\frac{B\displaystyle\int K(u) du+D}{K(u)},
\end{equation}
where $D$ is a constant of integration.

On the one hand, if $B \neq 0$, taking into account the definition of $A(u)$ given in \eqref{A}, we obtain
\begin{equation}
\frac{K(u)}{B\displaystyle\int K(u)\, du + D} = \frac{C'(u)}{C(u)} - \frac{K'(u)}{K(u)}.
\end{equation}
A subsequent integration then yields
\begin{equation}\label{CuCaso1}
C(u) = E\, K(u) \left( B \displaystyle\int K(u)\, du + D \right)^{1/B},
\end{equation}
where $E$ is a constant.

On the other hand, if $B = 0$, it follows from \eqref{Aint} that $A(u) = \frac{D}{K(u)}$ and 
using  \eqref{A}, we then have
\begin{equation}
\frac{K(u)}{D} = \frac{C'(u)}{C(u)} - \frac{K'(u)}{K(u)},
\end{equation}
which, upon integration, yields
\begin{equation}\label{CuCaso1B0}
C(u) = E\, K(u) \exp\left( \frac{1}{D} \displaystyle\int K(u)\, du \right),
\end{equation}
where $E$ is a constant.

In conclusion we have that  $C(u)$ must assume the form given in either \eqref{CuCaso1} or \eqref{CuCaso1B0}.

Consequently, the PDE \eqref{eq} admits a four-parameter Lie group with infinitesimal generators given in \eqref{infgen}, together with
\begin{equation}\label{Xbeta}
X_4 = x \pder[]{x} - 2 A(u) \pder[]{u},
\end{equation}
corresponding to the parameter $\beta$, where $A(u)$ is given by \eqref{Aint}.

\medskip

\noindent \textbf{Case 1.3:} $\rho \neq 0$, $\beta\neq 0$.\\
It should be noticed that, in order for \eqref{rhoKA} to hold, as $\rho \neq 0$, it is necessary that
\begin{equation}
4\big(K(u) A(u)\big)' = -K(u),
\end{equation}
then  we obtain
\begin{equation}\label{KAM}
K(u) A(u) = -\frac{1}{4} \displaystyle\int K(u)\, du + M,
\end{equation}
where $M$ is a constant of integration.  

Taking into account that $A$ is given by \eqref{A}, it then follows that
\begin{equation}\label{Cu-Caso13}
C(u) = \frac{N K(u)}{\left(4M - \displaystyle\int K(u)\, du \right)^4},
\end{equation}
where $N$ is a constant.

The form of $C(u)$ obtained above corresponds to a particular case of $C(u)$ given in \eqref{CuCaso1}, obtained by setting $M = D$, $N = 4^4 E$, and $B = -\frac{1}{4}$. This guarantees that equation \eqref{ux2-1} is satisfied.

Here, PDE \eqref{eq} admits a five- parameter group with infinitesimal generators \eqref{infgen}, \eqref{Xbeta} and
\begin{equation}\label{Xrho}
    X_5=x^2 \pder[]{x}-4x A(u) \pder[]{u}
\end{equation}
corresponding to the parameter $\rho$, where, in this case, $A$ is derived from \eqref{KAM} and it is given by
\begin{equation}
    A(u)=\frac{4M-\displaystyle\int K(u) du}{4K(u)}.
\end{equation}

In summary, the infinitesimal generators corresponding to this case are
\begin{equation}
\begin{cases}
X_1 = \frac{x}{2} \pder[]{x} + t \pder[]{t},\\[4pt]
X_2 = \pder[]{x}, \\[4pt]
X_3 = \pder[]{t},\\[4pt]
X_4 = x \pder[]{x} - 2 \frac{B\int K(u) du+D}{K(u)} \pder[]{u}, \\[4pt]
 X_5=x^2 \pder[]{x}-x \frac{4M-\int K(u) du}{K(u)} \pder[]{u}.
\end{cases}
\end{equation}

\begin{obs}
    In Case 1, the generators \(X_1\), \(X_2\), and \(X_3\) always exist, as long as the functions \(K(u)\) and \(C(u)\) satisfy $\frac{C(u)}{K(u)}$ no constant. 
    In addition, generator \(X_4\) appears when the specific relations \eqref{CuCaso1} or \eqref{CuCaso1B0} between \(K\) and \(C\) is satisfied. 
    Finally, generator \(X_5\) arises under the more restrictive condition \eqref{Cu-Caso13}.
\end{obs}

{\color{black}The five-dimensional Lie algebra generated by the infinitesimal generators is characterized by the commutation relations listed in Table \ref{tab:Commutator_table_Case1}. The Lie bracket is given by
$$[X_i,X_j] = X_i X_j - X_j X_i$$
with the generators interpreted as differential operators.}

\begin{table}[h!!]
\centering
\caption{Commutator table Case 1}
\label{tab:Commutator_table_Case1}
\begin{tabular}{|c|c|c|c|c|c|}
\hline
 - & \(X_1\) & \(X_2\) & \(X_3\) & \(X_4\) & \(X_5\)\\ \hline
\makecell{\(X_1\)} 
    & \makecell{\(0\)} & \makecell{\(-\frac{1}{2}X_2\)} & \makecell{\(-X_3\)} & \makecell{\(0\)} & \makecell{\(\frac{1}{2}X_5\)}\\ \hline

\makecell{\(X_2\)} 
    & \makecell{\(\frac{1}{2}X_2\)} & \makecell{\(0\)} & \makecell{\(0\)} & \makecell{\(X_2\)} & \makecell{\(2X_4\)}\\ \hline

\makecell{\(X_3\)} 
    & \makecell{\(X_3\)} & \makecell{\(0\)} & \makecell{\(0\)} & \makecell{\(0\)} & \makecell{\(0\)}\\ \hline

\makecell{\(X_4\)} 
    & \makecell{\(0\)} & \makecell{\(-X_2\)} & \makecell{\(0\)} & \makecell{\(0\)} & \makecell{\(X_5\)}\\ \hline

\makecell{\(X_5\)} 
    & \makecell{\(-\frac{1}{2}X_5\)} & \makecell{\(-2X_4\)} & \makecell{\(0\)} & \makecell{\(-X_5\)} & \makecell{\(0\)}\\ \hline
\end{tabular}
\end{table}

Then, for each infinitesimal generator we will obtain the corresponding one-parameter Lie group of transformations \eqref{groupLie} solving the following first order system of ODEs
\begin{equation}
\begin{cases}
\frac{dx^*}{d\epsilon}=\xi_1(x^*,t^*,u^*,\epsilon),\\
\frac{dt^*}{d\epsilon}=\xi_2(x^*,t^*,u^*,\epsilon),\\
\frac{du^*}{d\epsilon}=\eta(x^*,t^*,u^*,\epsilon),
\end{cases}
\end{equation}
with initial values $x^*=x$, $t^*=t$ and $u^*=u$ at $\epsilon=0$.

The resulting one-parameter groups associated with each infinitesimal generator are listed below.

\begin{itemize}
    \item $S_1:= \{x^*=x\, e^{\epsilon/2}, t^*=t\, e^{\epsilon}, u^*=u\}$
    \item $S_2 := \{x^*=x+\epsilon, t^*=t, u^*=u\}$
    \item $S_3 := \{x^*=x, t^*=t+\epsilon, u^*=u\}$
    \item $S_4 := \left\{x^*=x\, e^{\epsilon}, t^*=t, u^*=G^{-1}\!\left( G(u)\, e^{-2\epsilon}\right)\right\}$
    \item $S_5 := \left\{x^*=\dfrac{x}{1-x\epsilon}, t^*=t, u^*=H^{-1}\!\left( \tfrac{H(u)}{1-x\epsilon}\right)\right\}$
\end{itemize}
where $G(u)=\left(B\!\int K(u)\,du + D \right)^{1/B}$ and $H(u)=4M-\int K(u)\, du$ are assumed to be invertible.

\begin{obs}
    For the particular case  when \(C(u)\) is constant and \(K(u)\) remains variable, 
    the resulting symmetry algebra coincides with the one reported in \cite{BlKu1989}.
\end{obs}

\subsection*{Case 2:  $\tfrac{C(u)}{K(u)} = \alpha$ constant }

In this case, the PDE \eqref{eq} transforms into
\begin{equation}
    \alpha K(u)u_t = K'(u)u_x^2 + K(u)u_{xx},\label{eq_c2}
\end{equation}
and the equations \eqref{ti}, \eqref{ux}, \eqref{ux2}, \eqref{uxx} transform into

\begin{equation}
    \alpha\frac{\partial \eta}{\partial t} - \frac{\partial^2 \eta}{\partial x^2} = 0, \label{ti_c2}
\end{equation}
\begin{equation}
    2K'(u)\frac{\partial \eta}{\partial x} 
    + K(u)\left( \alpha\frac{\partial \xi_1}{\partial t} 
    + 2\frac{\partial^2 \eta}{\partial x \partial u} - \frac{\partial^2 \xi_1}{\partial x^2}\right) = 0, \label{ux_c2}
\end{equation}
\begin{equation}
    \eta\left( \frac{\left(K'(u)\right)^2}{K(u)} - K''(u) \right) 
    + K'(u)\left( 2\frac{\partial \xi_1}{\partial x} - \frac{\partial \eta}{\partial u} - \xi_2'\right) 
    - K(u)\frac{\partial^2 \eta}{\partial u^2} = 0, \label{ux2_c2}
\end{equation}
\begin{equation}
    2\frac{\partial \xi_1}{\partial x} - \xi_2' = 0. \label{uxx_c2}
\end{equation}

From \eqref{ux2_c2} and \eqref{uxx_c2}, we get

\begin{equation}
    \frac{\partial^2 \eta}{\partial u^2} = - \frac{\partial}{\partial u}\left(\eta \frac{K'(u)}{K(u)} \right),
\end{equation}
then, integrating with respect to $u$, $\eta$ must satisfy the following equation

\begin{equation}\label{n_c2}
    \eta = - \left[\frac{\partial\eta}{\partial u} + g(x,t)  \right]\frac{K(u)}{K'(u)}.
\end{equation}

From \eqref{uxx_c2}, we obtain

\begin{equation}\label{xi1_c2}
    \xi_1(x,t) = \frac{1}{2}\xi_2'(t)x + a(t)
\end{equation}
and taking into account \eqref{ux_c2}, \eqref{n_c2}, \eqref{xi1_c2} it follows

\begin{equation}\label{eqG}
    g(x,t) = \frac{\alpha}{8}\xi''_2 x^2 + \frac{\alpha}{2}a'(t)x + d(t).
\end{equation}

From \eqref{ti_c2} and \eqref{n_c2} we have that $\alpha g_t - g_{xx} = 0$ then by using \eqref{eqG} we obtain the following polynomial equation in $x$

\begin{equation*}
    \frac{\alpha^2}{8}\xi'''_2 x^2 + \frac{\alpha^2}{2}a''(t)x + \alpha d'(t) - \frac{\alpha}{4}\xi''_2 = 0.
\end{equation*}
We get $\xi'''_2(t) = 0$, $a''(t) = 0$ and $\alpha d'(t) - \frac{\alpha}{4}\xi''_2 = 0$. This leads to 

\begin{equation}
    \xi_2(t) = \mu t^2 + \gamma t + \nu,
\end{equation}
\begin{equation}
    a(t) = \zeta t + \theta,\quad d(t) = \frac{\mu}{2}t + \kappa.  
\end{equation}

Then, \eqref{eqG} remains

\begin{equation}
    g(x,t) = \frac{\alpha}{4}\mu x^2 + \frac{\alpha}{2}\zeta x + \frac{\mu}{2}t + \kappa,
\end{equation}
which substituting into \eqref{n_c2} leads to

\begin{equation}
    \eta = -\frac{\int K(u) du}{K(u)}\left(\frac{\alpha}{4}\mu x^2 + \frac{\alpha}{2}\zeta x + \frac{\mu}{2}t + \kappa \right).
\end{equation}

Consequently, the infinitesimals are given by
\begin{equation}
\left\{
\begin{array}{l}
\xi_1(x,t) = \mu tx + \frac{\omega}{2}x + \zeta t + \theta,\\[4pt]
\xi_2(t) = \mu t^2 + \omega t + \nu,\\[4pt]
\eta(x,t,u) = -\frac{\int K(u) du}{K(u)}\left(\frac{\alpha}{4}\mu x^2 + \frac{\alpha}{2}\zeta x + \frac{\mu}{2}t + \kappa\right).
\end{array}\right.
\end{equation}

In conclusion, for this case, the PDE \eqref{eq_c2} admits a six-parameter Lie group of transformations with infinitesimal generators, given by
\begin{equation}\label{infgen_c2}
\left\{
\begin{array}{l}
\overline{X}_1 = tx\pder[]{x} + t^2\pder[]{t} -\left(\frac{\alpha x^2}{4} + \frac{t}{2}\right) \frac{\int K(u) du}{K(u)} \pder[]{u}, \\[6pt]
\overline{X}_2 = \frac{x}{2}\pder[]{x} + t\pder[]{t}, \\[6pt]
\overline{X}_3 = t\pder[]{x} - \frac{\alpha x\int K(u) du}{2K(u)} \pder[]{u}, \\[6pt]
\overline{X}_4 = \pder[]{x}, \\[6pt]
\overline{X}_5 = \pder[]{t}, \\[6pt]
\overline{X}_6 = -\frac{\int K(u) du}{K(u)} \pder[]{u},
\end{array}\right.
\end{equation}
corresponding to the parameters $\mu$, $\omega$, $\zeta$, $\theta$, $\nu$ and $\kappa$ respectively.\\

{\color{black}The commutator table for the Lie algebra arising from the infinitesimal generators \eqref{infgen_c2} is presented in Table \ref{tab:Commutator_table_Case2}.}
\begin{table}[h!!]
\centering
\renewcommand{\arraystretch}{1.25}
\caption{Commutator table Case 2}
\label{tab:Commutator_table_Case2}
\begin{tabular}{|c|c|c|c|c|c|c|}
\hline 
$-$ & $\overline{X}_1$ & $\overline{X}_2$ & $\overline{X}_3$ & $\overline{X}_4$ & $\overline{X}_5$ & $\overline{X}_6$ \\ \hline

$\overline{X}_1$ 
& $0$ & $-\overline{X}_1$ & $0$ & $-\overline{X}_3$ & $-2\overline{X}_2-\frac{1}{2}\overline{X}_6$ & $0$ \\ \hline

$\overline{X}_2$ 
& $\overline{X}_1$ & $0$ & $\frac{1}{2}\overline{X}_3$ & $-\frac{1}{2}\overline{X}_4$ & $-\overline{X}_5$ & $0$ \\ \hline

$\overline{X}_3$ 
& $0$ & $-\frac{1}{2}\overline{X}_3$ & $0$ & $-\frac{\alpha}{2}\overline{X}_6$ & $-\overline{X}_4$ & $0$ \\ \hline

$\overline{X}_4$ 
& $\overline{X}_3$ & $\frac{1}{2}\overline{X}_4$ & $\frac{\alpha}{2}\overline{X}_6$ & $0$ & $0$ & $0$ \\ \hline

$\overline{X}_5$ 
& $2\overline{X}_2+\frac{1}{2}\overline{X}_6$ & $\overline{X}_5$ & $\overline{X}_4$ & $0$ & $0$ & $0$ \\ \hline

$\overline{X}_6$ 
& $0$ & $0$ & $0$ & $0$ & $0$ & $0$ \\ \hline

\end{tabular}
\end{table}

The corresponding one-parameter Lie groups of transformations generated by each infinitesimal operator are presented below
\begin{itemize}
    \item $\overline{S}_1:= \left\{x^*=\frac{x}{1-\epsilon t}, t^*=\frac{t}{1-\epsilon t}, u^*=I^{-1}\left(I(u)\sqrt{1-\epsilon t}\exp\left(-\frac{\alpha\epsilon x^2}{4(1-\epsilon t)}\right)\right)\right\}$
    \item $\overline{S}_2 := \{x^*=x\, e^{\epsilon/2}, t^*=t\, e^{\epsilon}, u^*=u\}$
    \item $\overline{S}_3 := \left\{x^*=\epsilon t + x, t^*=t, u^*=I^{-1}\left(I(u)\exp\left(\frac{-\alpha\epsilon^2t}{4}-\frac{\alpha\epsilon x}{2}\right)\right)\right\}$
    \item $\overline{S}_4 := \{x^*=x+\epsilon, t^*=t, u^*=u\}$
    \item $\overline{S}_5 := \{x^*=x, t^*=t+\epsilon, u^*=u\}$
    \item $\overline{S}_6 := \{x^*=x, t^*=t, u^*=I^{-1}\left(I(u)\exp\left(\epsilon\right)\right)\}$
\end{itemize}

where $I(u)= \!\int K(u)\,du$ is assumed to be invertible.

\begin{obs}
    Note that the one-parameter groups generated by $\overline{X}_2$, $\overline{X}_4$, and $\overline{X}_5$ coincide with those of $X_1$, $X_2$, and $X_3$ from Case 1, respectively.
\end{obs}

\section{Invariant solutions \(u = \Theta(x,t)\)  }
In this section, we determine the invariant solutions $u=\Theta(x,t)$ of PDE~\eqref{eq} corresponding to each of its admitted infinitesimal generators.

Recall that a function $u=\Theta(x,t)$ is an invariant solution associated with an infinitesimal generator $X$ if and only if
\begin{equation}\label{cond-inv}
    X\big(u-\Theta(x,t)\big)=0 \qquad \text{whenever} \qquad u=\Theta(x,t).
\end{equation}
Using this criterion, we now determine the invariant solutions of PDE \eqref{eq} corresponding to the infinitesimal generators $X_i$, $i=1,2,3,4,5$, for Case~1, and to $\overline{X}_i$, $i=1,2,3,4,5,6$, for Case~2.

\subsection{Case 1}
\noindent \textbf{Generator \(X_1\).}
We begin with the generator \(X_1=\frac{x}{2}\pder{x}+t\pder{t}\), given by \eqref{infgen}.
If \(u = \Theta(x,t)\)  is an invariant solution to \eqref{eq}, then from \eqref{cond-inv} it must satisfy

\begin{equation}
    \frac{x}{2}\pder[\Theta]{x} + t\pder[\Theta]{t} = 0.
\end{equation}
The corresponding characteristic equations show that \(u\) remains constant along the characteristic curves and that
\begin{equation}
    \frac{dx}{x/2} = \frac{dt}{t},
\end{equation}
their solution gives the similarity variable
\begin{equation}
    \xi = \frac{x}{\sqrt{t}},
\end{equation}
hence, an invariant solution $u$ of \eqref{eq} can be written as \(u = \phi_1(\xi)\). Substituting this expression into the PDE~\eqref{eq} yields the reduced ordinary differential equation satisfied by $\phi_1$:
\begin{equation}\label{eq_dif_phi_delta}
    2K(\phi_1)\,\phi_1''
    + 2K'(\phi_1)\,(\phi_1')^2
    + \xi\,C(\phi_1)\,\phi_1' = 0,
\end{equation}
which can be written equivalently as
\begin{equation}
    2\big(K(\phi_1)\,\phi_1'\big)' + \xi\,C(\phi_1)\,\phi_1' = 0.
\end{equation}

Setting $v = K(\phi_1)\,\phi_1'$, we obtain
\begin{equation}
    v' + \frac{\xi}{2}\,\frac{C(\phi_1)}{K(\phi_1)}\,v = 0.
\end{equation}
Solving this ordinary differential equation gives
\begin{equation}
    v(\xi)=\widetilde{D}\,
    \exp\!\left(
        -\frac{1}{2}\int \xi\,\frac{C(\phi_1(\xi))}{K(\phi_1(\xi))}\,d\xi
    \right),
\end{equation}
where $\widetilde{D}$ is a constant.

Therefore,
\begin{equation}\label{phi1-eqint}
    \phi_1(\xi)
    =
    \widetilde{E}
    + \widetilde{D}
      \int
      \frac{1}{K(\phi_1(\xi))}
      \exp\!\left(
        -\frac{1}{2}\int 
        \xi\,\frac{C(\phi_1(\xi))}{K(\phi_1(\xi))}\,d\xi
      \right)\,d\xi,
\end{equation}
where $\widetilde{E}$ is a constant, which shows that $\phi_1$ must satisfy the above integral equation.\\
\ \\
\textbf{Generator \(X_2\).}
Next, we consider the generator \(X_2=\frac{\partial}{\partial x}\), given by \eqref{infgen}. In this case, the invariance condition yields
\begin{equation}
   \pder[\Theta]{x} = 0,
\end{equation}
hence, \(u = \phi_2(t)\) and from \eqref{eq} it follows that \(u \equiv u_0\), constant.\\
\ \\
\textbf{Generator \(X_3\).}
For the generator $X_3=\frac{\partial}{\partial t} $ given by \eqref{infgen}, the invariant solution $u=\Theta(x,t)$ of \eqref{eq}, must satisfy 
\begin{equation}
   \pder[\Theta]{t} = 0,
\end{equation}
that is, \(u = \phi_3(x)\). Substituting this expression into \eqref{eq}, we obtain
\[
\left(K(\phi_3(x))\phi'_3(x)\right)' = 0.
\]
Consequently $\phi_3$ satisfies the following integral equation,
\begin{equation}\label{caso2-theta}
   \phi_3(x) = \int \frac{u_1}{K(\phi_3(x))}\,dx,
\end{equation}
where $u_1$ is constant.\\
\ \\
\textbf{Generator \(X_4\).}
Next, \(X_4=x\pder{x}-2\frac{B\int K(u)du+D}{K(u)}\pder{u}\) with the invariant solutions \(u = \Theta(x,t)\) that must satisfy the invariance condition \(X_4(u - \Theta) = 0\), that is,
\begin{equation}
    x\pder[\Theta]{x} + \frac{2}{K(\Theta)}\left(B\int K(u)\,du + D\right) = 0.
\end{equation}
The associated characteristic equations imply that \(\xi=t\) is the similarity variable and
\begin{equation}\label{chxbeta}
    \frac{2\,dx}{-x} = \frac{K(u)\,du}{B\int K(u)\,du + D}.
\end{equation}
If \(B \neq 0\), we obtain the invariant relation
\begin{equation}\label{InvariantXbeta-Caso1}
    x\,\phi_4(t) = \left(B\int K(u)\,du + D\right)^{-\tfrac{1}{2B}}.
\end{equation}
Differentiating \eqref{InvariantXbeta-Caso1} with respect to \(x\) gives
\begin{equation}\label{phi}
    \phi_4(t) = \frac{1}{2}\left(B\int K(u)\,du + D\right)^{-\tfrac{1}{2B}-1} K(u)\,\pder[u]{x},
\end{equation}
and differentiating once more with respect to $x$ yields
\begin{equation}\label{phi2}
\begin{aligned}
    0 = &\, \tfrac{B}{2}\!\left(-\tfrac{1}{2B} - 1\right)
    \left(B\int K(u)\,du + D\right)^{-\tfrac{1}{2B}-2} K^2(u)\,\left( \pder[u]{x}\right)^2 \\
    &+ \tfrac{1}{2}\left(B\int K(u)\,du + D\right)^{-\tfrac{1}{2B}-1}
    \Big(K'(u)\,\left( \pder[u]{x}\right)^2 + K(u)\,\pder[^2u]{x^2}\Big).
\end{aligned}
\end{equation}
Similarly, differentiating \eqref{InvariantXbeta-Caso1} with respect to \(t\) gives
\begin{equation}\label{xphi'}
    x\,\phi_4'(t) = \frac{1}{2}\left(B\int K(u)\,du + D\right)^{-\tfrac{1}{2B}-1} K(u)\,\pder[u]{t}.
\end{equation}
On one hand, from \eqref{InvariantXbeta-Caso1} and \eqref{xphi'} we obtain that
\begin{equation}\label{eq:Kut}
\begin{aligned}
    K(u)\,u_t
    &= \frac{2x\,\phi_4'(t)}
        {\left(B\!\int K(u)\,du + D\right)^{-\tfrac{1}{2B}-1}} = \frac{2\,\phi_4'(t)}{\phi_4(t)}
        \left(B\!\int K(u)\,du + D\right).
\end{aligned}
\end{equation}
On the other hand,  taking into account \eqref{phi} and \eqref{phi2},  we get that
\begin{equation}\label{eq:Kuxx_phi}
\begin{aligned}
    K'(u)\!\left(\pder[u]{x}\right)^{\!2} + K(u)\,\pder[^2u]{x^2}
    &= (2B + 1)\, 2\, \phi_4^2(t)\,
       \left(B\!\int K(u)\,du + D\right)^{1 + \tfrac{1}{B}}.
\end{aligned}
\end{equation}
Taking PDE~\eqref{eq} into account, together with the expression for \(C\) given in \eqref{CuCaso1}, it follows that
\begin{equation}
    E \frac{\phi_4'(t)}{\phi_4^3(t)}=2B+1.
\end{equation}
Integrating with respect to $t$ gives
\begin{equation}\label{Phicuad}
    \phi_4^2(t)=\frac{E}{QE-(2+4B)t},
\end{equation}
where $Q$ is a constant of integration.

In summary, $u$ is implicitly defined by \eqref{InvariantXbeta-Caso1}, with $\phi_4$ given by \eqref{Phicuad}.

In the case that $B=0$  from \eqref{chxbeta} we get the invariant relation
\begin{equation}\label{u-xbeta-B=0}
    x \phi_4(t)= \exp \left( -\tfrac{1}{2D}\displaystyle\int K(u) du \right) .
\end{equation}
Upon differentiating with respect to \(x\), twice with respect to \(x\), and with respect to \(t\), we obtain, respectively,
\begin{equation}\label{deriv-xbeta}
\begin{aligned}
&\phi_4(t) 
    = \frac{-1}{2D} \exp\Biggl(-\frac{1}{2D}\int K(u)\,du\Biggr)\, K(u)\,\pder[u]{x}, \\[2mm]
&0 
    = \exp\Biggl(-\frac{1}{2D}\int K(u)\,du\Biggr)
       \Biggl[
       \frac{1}{4D^2} K^2(u) \left(\pder[u]{x}\right)^2
        -\frac{1}{2D} \left(K'(u) \left(\pder[u]{x}\right)^2 + K(u)\,\pder[^2 u]{x^2}\right)
       \Biggr], \\[1mm]
&\phi_4'(t)\,x 
    = -\frac{1}{2D} \exp\Biggl(-\frac{1}{2D}\int K(u)\,du\Biggr)\, K(u)\, \pder[u]{t}.
\end{aligned}
\end{equation}
As a consequence we get
\begin{equation}\label{eq:Kut_phi_multiline}
\begin{aligned}
K(u)\,\pder[u]{t} 
    &= - 2D\,\phi_4'(t)\,x \,\exp\Biggl(\frac{1}{2D}\int K(u)\,du\Biggr) =- \frac{2D\,\phi_4'(t)}{\phi_4(t)}  .
\end{aligned}
\end{equation}
In a similar manner, taking into account \eqref{deriv-xbeta} we get
\begin{equation}\label{eq:Kuxx_phi2}
K'(u)\,\left(\pder[u]{x}\right)^2 + K(u)\,\pder[^2 u]{x^2} 
    = \frac{1}{2D} \, K^2(u) \,\left(\pder[u]{x}\right)^2 
    = 2D\, \phi_4^2(t) \,\exp\Biggl(\frac{1}{D}\int K(u)\,du\Biggr).
\end{equation}
From PDE \eqref{eq} and the fact that $C$ is given by \eqref{CuCaso1B0},  it follows that
\begin{equation}
   E \frac{\phi_4'(t)}{\phi_4^3(t)}=-1.
\end{equation}
Integrating with respect to $t$ gives
\begin{equation}\label{Phicuad0}
    \phi_4^2(t)=\frac{E}{2t+Q E},
\end{equation}
where $Q$ is a constant of integration.

Therefore we have the solution $u$  implicitly defined by \eqref{u-xbeta-B=0}, where $\phi_4$ is given by \eqref{Phicuad0}.\\
\ \\
\textbf{Generator \(X_5\).}
Finally, for \(X_5=x^2\pder{x}-x\frac{4M-\int K(u)du}{K(u)}\pder{u}\) given by \eqref{Xrho} the invariant solution \(u = \Theta(x,t)\) of PDE~\eqref{eq} corresponding to \(X_5\) satisfy \(X_5(u - \Theta) = 0\) when \(u = \Theta\), which gives
\begin{equation}\label{chXrho}
    x^2 \pder[\Theta]{x} = x \left( \frac{1}{K(u)} \int K(u)\,du - 4 \frac{M}{K(u)} \right).
\end{equation}
The associated characteristic equations imply that \(t\) is constant and
\begin{equation}
    \frac{dx}{x} = \frac{K(u)\,du}{ \left( \int K(u)\,du - 4M \right)}.
\end{equation}
Then we get the invariant form
\begin{equation}
    x= \left( \displaystyle\int K(u)\,du - 4M \right)  \phi_5(t).
\end{equation}
Upon differentiating with respect to \(x\), twice with respect to \(x\), and with respect to \(t\), we obtain, respectively,
\begin{equation}\label{eq:three_relations_Xrho}
\begin{aligned}[t]
&  K(u)\, \pder[u]{x} \phi_5(t)=1, \\[1mm]
& \left(K'(u) \left(\pder[u]{x} \right)^2 + K(u)\, \pder[^2 u]{x^2}\right) \phi_5(t)=0, \\[1mm]
 &  K(u)\, \pder[u]{t}\, \phi_5(t) + \left( \int K(u)\,du - 4M \right) \phi_5'(t)=0.
\end{aligned}
\end{equation}
From PDE~\eqref{eq} it follows that $\phi_5'(t)=0$ and therefore
\begin{equation}
    \phi_5(t) = u_2,
\end{equation}
where \(u_2\) is a constant.

In Table 3 we {\color{black}summarize} the invariant solutions of equation \eqref{eq} for Case 1  with the corresponding $K$ and $C$.

\begin{table}[h!]
\centering
\caption{Invariant solutions for Case 1}
\label{tab:invariant_solutions}
\begin{tabular}{|c|c|}
\hline
Generator & Solution \(u\) \\ \hline

\makecell{\(X_1 = \frac{x}{2}\pder[]{x} + t \pder[]{t}\)} 
&
\makecell{ \\ \(u = \phi_1(\xi)\) ,\quad \(\xi = x/\sqrt{t}\)\\ where\\
\(\phi_1(\xi)= \widetilde{E}+\widetilde{D}\displaystyle\int\tfrac{\exp\!\left(-\frac{1}{2}\int \xi\,\tfrac{C(\phi_1(\xi))}{K(\phi_1(\xi))}\,d\xi\right)}{K(\phi_1(\xi))}
\,d\xi\)\\$ $}
\\  \hline
&\\
\makecell{\(X_2 = \pder[]{x} \)} 
&
\makecell{ \(u=u_0\) constant \\ $ $}
\\ \hline
&\\
\makecell{\(X_3 = \pder[]{t} \)} 
&
\makecell{\(u=\phi_3(x)\) \\ 
where \\  
\(\phi_3(x)=\int \frac{u_1}{K(\phi_3(x))}  dx, \quad u_1\) constant \\ $ $}
\\ \hline
&\\
\makecell{\(X_4 = x\pder[]{x} - \frac{2}{K(u)} (B\int K(u)\,du + D)\pder[]{u}\) \\
where\\
$C(u)=E\,K(u)\left(B \int K(u)\,du + D\right)^{\tfrac{1}{B}}$ \\
if $B \neq 0$; \, \, or\\ \\
$C(u)=E\,K(u)\exp\!\left(\tfrac{1}{D}\int K(u)\,du\right),
$\\
if $B=0$ \\ $ $}
&
\makecell{
If $B \neq 0$, \color{black}{$u$ is given implicitly by}\\
$\phi_4(t)\, x = \left(B\!\int K(u)\,du + D\right)^{-1/(2B)}$ \\[4pt]
where $\displaystyle 
\phi_4^2(t)=\tfrac{E}{Q E - (2+4B)t}$ \\[6pt]
If $B = 0$, \color{black}{$u$ is given implicitly by}\\
$x \phi_4(t)= \exp \left( \tfrac{-1}{2D}\int K(u) du \right)$ \\
where $\phi_4^2(t)=\tfrac{E}{Q E +2t}$ 
\\ $ $}
 \\ \hline 
& \\
\makecell{$X_5=x^2 \pder[]{x}-\frac{x}{K} (4M-\int K(u)du) \pder[]{u}$\\
where\\
$C(u)=\frac{N\,K(u)}{\left(4M-\int K(u)\,du\right)^4}$ \\ $ $}
&
\makecell{ $u$ is given implicitly by\\$x= \left( \int K(u)\,du - 4M \right)  u_2$ \\ 
with \(u_2\) constant \\ $ $}
 \\ \hline 

\end{tabular}
\end{table}
\newpage

\subsection{Case 2}

\noindent\textbf{Generator \(\overline{X}_1\).}
Proceeding now with the generator \(\overline{X}_1\) given by \eqref{infgen_c2}, if \(u = \Theta(x,t)\)  is an invariant solution to \eqref{eq}, then it must satisfy $\overline{X}_1(u-\Theta(x,t))=0$ which implies
\begin{equation}
    xt\pder[\Theta]{x} + t^2\pder[\Theta]{t} + \left(\frac{\alpha x^2}{4} + \frac{t}{2}\right) \frac{\int K(u) du}{K(u)} = 0.
\end{equation}

The corresponding characteristics equations show that
\begin{equation}
    \frac{dx}{xt} = \frac{dt}{t^2} = -\frac{K(u)du}{\left(\frac{\alpha x^2}{4} + \frac{t}{2}\right)\int K(u) du} .
\end{equation}
On one hand
\begin{equation}
    \frac{dx}{xt} = \frac{dt}{t^2}
\end{equation}
gives the similarity variable
\begin{equation}
    \xi = \frac{x}{t},
\end{equation}
while from
\begin{equation}
    \frac{dt}{t^2} = -\frac{K(u)du}{\left(\frac{\alpha x^2}{4} + \frac{t}{2}\right) \int K(u) du},
\end{equation}
we obtain the relation
\begin{equation}\label{InvariantXmu-Caso2}
    \int K(u) du = \frac{1}{\sqrt{t}}\psi_1(\xi)\exp\left(-\frac{\alpha\xi^2t}{4}\right).
\end{equation}
Differentiating \eqref{InvariantXmu-Caso2} with respect to \(x\) gives
\begin{equation}\label{kux}
    K(u)\pder[u]{x} = \frac{1}{\sqrt{t}}\exp\left(-\frac{\alpha\xi^2t}{4}\right)\left[\frac{1}{t}\psi_1'(\xi) - \frac{\alpha\xi}{2}\psi_1(\xi)\right],
\end{equation}
and differentiating once more with respect to $x$ yields 
\begin{equation}\label{k'ux2+kuxx}
    \frac{\partial }{\partial x}\left(K(u)\pder[u]{x}\right) = \frac{1}{\sqrt{t}}\exp\left(-\frac{\alpha\xi^2t}{4}\right)\left[\frac{1}{t^2}\psi_1''(\xi)- \frac{\alpha\xi}{t}\psi_1'(\xi) + \left(\frac{\alpha^2\xi^2}{4}-\frac{\alpha}{2t}\right)\psi_1(\xi)\right].
    \end{equation}

Similarly, differentiating \eqref{InvariantXmu-Caso2} with respect to \(t\) gives
\begin{equation}\label{kut}
    K(u)\pder[u]{t} = \frac{1}{\sqrt{t}}\exp\left(-\frac{\alpha\xi^2t}{4}\right)\left[\frac{\alpha\xi^2}{4}\psi_1(\xi) - \frac{1}{2t}\psi_1(\xi) - \frac{\xi}{t}\psi_1'(\xi)\right].
\end{equation}
Substituting \eqref{k'ux2+kuxx} and \eqref{kut} into PDE \eqref{eq_c2} yields $\psi_1''(\xi) = 0$, then we find that solution $u$ is implicitly defined by
\begin{equation}\label{psi1_eq}
    \int K(u) du = \frac{a\xi +b}{\sqrt{t}}\exp\left(-\frac{\alpha\xi^2t}{4}\right), \quad a,b\in \mathbb{R}.
\end{equation}
\ \\
\textbf{Generator \(\overline{X}_2\).}
By analogy with the invariant solution \(u=\phi_1(\xi)\), where 
\(\xi = x/ \sqrt{t}\) associated with the generator \(X_1\), we obtain the invariant solution  \(u=\psi_2(\xi)\) for the infinitesimal generator $\overline{X}_2=\frac{x}{2}\pder{x}+t\pder{t}$ given by \eqref{infgen_c2}. We have \(\psi_2\) satisfies

\begin{equation}
    2K(\psi_2(\xi))\,\psi_2''(\xi)
    + 2K'(\psi_2(\xi))\,\big(\psi_2'(\xi)\big)^2
    + \xi\alpha K(\psi_2(\xi))\,\psi_2'(\xi) = 0,
\end{equation}
therefore, $\psi_2(\xi)$ must satisfy the integral equation
\begin{equation}\label{psi2}
    \psi_2(\xi) = \widetilde{E} + \widetilde{D}\int \frac{1}{K(\psi_2(\xi))}\exp\left(-\frac{\alpha\xi^2}{4}\right)d\xi.
\end{equation}
\ \\
\textbf{Generator \(\overline{X}_3\).}
For the generator \(\overline{X}_3=t\pder{t}-\frac{\alpha\int K(u)du}{2K(u)}\pder{u}\) given by \eqref{infgen_c2}, if \(u = \Theta(x,t)\)  is an invariant solution to \eqref{eq_c2}, then it must satisfy

\begin{equation}
    t\pder[\Theta]{x} + \frac{\alpha x}{2}\frac{\int K(u) du}{K(u)} = 0,
\end{equation}
then the associated characteristic equations imply that \(\xi=t\) is the similarity variable and
\begin{equation}
    \frac{dx}{t} = -\frac{2K(u)du}{\alpha x\int K(u) du}.
\end{equation}
We obtain the relation
\begin{equation}\label{InvariantXzeta-Caso2}
    \int K(u)du = \psi_3(t)\exp\left(-\frac{\alpha x^2}{4t}\right),
\end{equation}
for some function $\psi_3(t)$.

Differentiating \eqref{InvariantXzeta-Caso2} with respect to \(x\) gives

\begin{equation}\label{kux_z}
    K(u)\pder[u]{x} = -\frac{\alpha x}{2t}\psi_3(t)\exp\left(-\frac{\alpha x^2}{4t}\right)
\end{equation}
and differentiating once more with respect to $x$ yields
\begin{equation}\label{k'ux2+kuxx_z}
    K'(u)\left(\pder[u]{x}\right)^2 + K(u)\pder[^2u]{x^2} = \psi_3(t)\exp\left(-\frac{\alpha x^2}{4t^2}\right)\left[\frac{\alpha^2x^2}{4t^2}-\frac{\alpha}{2t}\right].
\end{equation}
Similarly, differentiating \eqref{InvariantXzeta-Caso2} with respect to \(t\) gives
\begin{equation}\label{kut_z}
    K(u)\pder[u]{t} = \exp\left(-\frac{\alpha x^2}{4t}\right)\left[\psi_3'(t) + \frac{\alpha x^2}{4t^2}\psi_3(t)\right].
\end{equation}
Substituting \eqref{k'ux2+kuxx_z} and \eqref{kut_z} into the PDE \eqref{eq_c2} yields 
\begin{equation}
    \psi_3' + \frac{\psi_3}{2t} = 0
\end{equation}
which implies that $\psi_3(t) = \frac{a}{\sqrt{t}},\text{ } a \in \mathbb{R}$, and we obtain that the solution $u$ is implicitly defined by
\begin{equation}\label{u-X1raya}
    \int K(u) du = \frac{a}{\sqrt{t}}\exp\left(-\frac{\alpha x^2}{4t}\right).
\end{equation}
\ \\
\textbf{Generator \(\overline{X}_4\).}
Next, considering \(\overline{X}_4=\pder{x}\) given by \eqref{infgen_c2}. The invariant solution \(u = \Theta(x,t)\) satisfies
\begin{equation}
    \frac{\partial\Theta}{\partial x}=0,
\end{equation}
this implies that $u=\psi_4(t)$ and 
\begin{equation}
    \alpha K(u)\psi_4'(t)=0.
\end{equation}
Hence, the solution is that $u(x,t) = \overline{u}_0$ constant.

\medskip

\noindent\textbf{Generator \(\overline{X}_5\).}
For \(\overline{X}_5=\pder{t}\) given by \eqref{infgen_c2} we have
\begin{equation}
     \frac{\partial\Theta}{\partial t}=0,
\end{equation}
this implies that $u=\psi_5(x)$. Using this in the PDE \eqref{eq_c2}
\begin{equation}
    0= K'(\psi_5(x))(\psi_5'(x))^2+K(\psi_5(x))\psi_5''(x),
\end{equation}
then $\psi_5$ must satisfy the integral equation
\begin{equation}\label{psi5}
    \psi_5(x) = a\int\frac{dx}{K(\psi_5(x))}+b, \quad a,b \in \mathbb{R}.
\end{equation}
\ \\
\textbf{Generator \(\overline{X}_6\).}
Finally, for generator \(\overline{X}_6=-\frac{\int K(u)du}{K(u)}\pder{u}\) given by \eqref{infgen_c2}, if \(u = \Theta(x,t)\)  is an invariant solution to \eqref{eq_c2}, then it must satisfy
\begin{equation}
     \frac{\int K(u)du}{K(u)}=0,
\end{equation}
which is not satisfied for any non zero $K$.

In Table 4 we present the invariant solutions of (PDE) \eqref{eq_c2}.
\begin{table}[h!]
\centering
\caption{Invariant solutions for Case 2}
\label{tab:invariant_solutions_case2}
\begin{tabular}{|c|c|}
\hline
Generator & Solution \(u\) \\ \hline
\makecell{\(\overline{X}_1 = tx\pder[]{x} + t^2\pder[]{t} -\left(\frac{\alpha x^2}{4} + \frac{t}{2}\right) \frac{\int K(u) du}{K(u)} \pder[]{u}\)} 
    & \makecell{\\ $u$ is given implicitly by \\\(\int K(u) du = \frac{1}{\sqrt{t}} \psi_1(\xi)\exp\left(-\frac{\alpha \xi^2 t}{4}\right)\)\\ where $\psi_1(\xi)= a\xi + b$, $\xi=x/t $ \\ $ $ } \\  \hline

\makecell{ \(\overline{X}_2 = \frac{x}{2}\pder[]{x} + t \pder[]{t}\)}
    & \makecell{\\ \(u = \psi_2(\xi)\),  \(\xi = x/\sqrt{t}\) \\ satisfies \\ $\psi_2(\xi) = \widetilde{E} +\widetilde{D}\int \frac{1}{K(\psi_2(\xi))}\exp\left(-\frac{\alpha\xi^2}{4}\right)d\xi$ \\ $ $} \\ \hline

\makecell{\(\overline{X}_3 = t\pder[]{x} - \frac{\alpha x\int K(u) du}{2K(u)} \pder[]{u}\) }
    & \makecell{\\ $u$ is given implicitly by \\ \(\int K(u)du=\frac{a}{\sqrt{t}}\exp\left(-\frac{\alpha x^2}{4t}\right)\)\\ $ $} \\ \hline

\makecell{\(\overline{X}_4 = \pder[]{x}\)} 
    & \makecell{\\ \(u\equiv \overline{u}_0\) \\ $ $} \\ \hline

\makecell{\(\overline{X}_5 = \pder[]{t}\)} 
    & \makecell{\\ \(u=\psi_5(x)\) satisfies\\\(\psi_5(x) = a\int\frac{dx}{K(\psi_5(x))}+b\) \\ $ $} \\ \hline

\makecell{\\ \(\overline{X}_6 = -\frac{\int K(u)du}{K(u)}\pder[]{u}\) \\ $ $}
    & \makecell{-} \\ \hline

\end{tabular}
\end{table}

\section{Particular differential equations }
Here we present three particular partial differential equations of the type given in \eqref{eq} with variable coefficients found in the literature. For such equations, we will seek invariant solutions using the results obtained in Section 3.

\subsection{{\color{black}Power-type $C(u)$ and constant $K(u)$}}
First we consider the nonlinear PDE \eqref{eq} given by
\begin{equation}\label{ejemplo1}
    \frac{u_t}{u^2}= k u_{xx}
\end{equation}
where $C(u)=\frac{1}{u^2}$ and $K(u)=k$. 
One--phase and two--phase Stefan problems governed by this equation were considered in
\cite{DeSa2002,DeSa2004}.
In this case, $C$ and $K$ satisfy relation \eqref{CuCaso1} for $B=-\frac{1}{2}$, $E=\frac{k}{4}$ and $D=0$. Based on the results obtained in Sections 2 and 3 for the case 1, we have that equation \eqref{ejemplo1} is admitted by the generators $X_1$, $X_2$, $X_3$ and $X_4$ given in Table \ref{tab:invariant_solutions} and the corresponding invariant solutions are:
\begin{itemize}
    \item $u=\phi_1(\xi)$  where $\phi_1$ must be a solution of the integral equation \eqref{phi1-eqint}, i.e.,
    \begin{equation}
    \phi_1(\xi)
    =
    \widetilde{E}
    + \frac{\widetilde{D}}{k}
      \int
      \exp\!\left(
        -\frac{1}{2k}\int 
        \frac{\xi}{\phi_1^2(\xi)}\,d\xi
      \right)\,d\xi,
\end{equation}
where $\xi=\frac{x}{\sqrt{t}}$.

\item $u\equiv u_0$ constant.

\item $u=\phi_3(x)$ where $\phi_3 $ satisfies equation \eqref{caso2-theta}, that is
\begin{equation}
    \phi_3(x)= \frac{u_1}{k}x+\bar{u}_1, \quad u_1, \bar{u}_1 \in \mathbb{R}.
\end{equation}

\item $u$ satisfies equation  \eqref{InvariantXbeta-Caso1}, that is
\begin{equation}
    u=-2\frac{x}{k}\phi_4(t),
\end{equation}
where $\phi_4(t)$ given by  \eqref{Phicuad} becomes constant, i.e. 
\begin{equation}
u=-2\frac{x}{k}\bar{u}_2,\qquad \bar{u}_2\in\mathbb{R}.\end{equation}
\end{itemize}

\subsection{{\color{black}Storm's condition for $C(u)$ and $K(u)$}}

Now, we consider equation \eqref{eq} in the case where the functions \(C\) and \(K\) satisfy the Storm’s condition
\cite{BrNa2014-1,Br2018,KnPh1974,NaTa2000,St1951}:
\begin{equation}\label{storm}
\frac{d}{du}\,\sqrt{\frac{C(u)}{K(u)}}
= \lambda\,C(u) , \quad \lambda>0,
\end{equation}
which is equivalent to \begin{equation}
   C(u)=\frac{K(u)}{4\lambda^{2}}\left(-\frac{1}{2}\int K(u)du + D\right)^{-2},
\end{equation} with $D$ constant.
This condition implies that \(C\) and \(K\) are related through \eqref{CuCaso1},
with \(E=\frac{1}{4\lambda^2}\) and \(B=-\frac{1}{2}\).
In view of the results established in Sections~2 and 3 for case 1, equation \eqref{eq} admits the symmetry generators
\(X_1\), \(X_2\), \(X_3\), and \(X_4\), whose corresponding invariant solutions are given by
\begin{itemize}
    \item $u=\phi_1(\xi)$  where $\phi_1$ must be a solution of the integral equation \eqref{phi1-eqint}, i.e.,
    \begin{equation}\label{fi1}
    \phi_1(\xi)
    =
    \widetilde{E}
    + \widetilde{D}
      \int
      \tfrac{1}{K(\phi_1(\xi))}
      \exp\!\left(
        -\tfrac{1}{8\lambda^2}\int 
        \tfrac{ \xi\, }{\left(D-\tfrac{1}{2}\int K(\phi_1(\xi))d\xi  \right)^2} \,d\xi
      \right)\,d\xi,
\end{equation}
where $\xi=\frac{x}{\sqrt{t}}$.

\item $u\equiv u_0$ constant.

\item $u=\phi_3(x)$ where $\phi_3 $ satisfies the integral equation \begin{equation}\label{fi3}
    \phi_3(x) = \int \frac{u_1}{K(\phi_3(x))}\,dx.
\end{equation}

\item $u$ is defined implicitly by \eqref{InvariantXbeta-Caso1} where $\phi_4$  given by  \eqref{Phicuad} turns out to be the  constant $\frac{1}{\sqrt{Q}}$, that is
\begin{equation}\label{fi4}
    \int K(u) du= \frac{-2x}{\sqrt{Q}} +2D.
\end{equation}
\end{itemize}

{\color{black}For simple metals \cite{St1951}, sufficient conditions for the validity of \eqref{storm} are $KC$ constant and the exponential behavior of $K$ and $C$ that is to say
\begin{equation}
K(u)=k_0\exp(-\widetilde{A}u),\qquad C(u)=c_0\exp(\widetilde{A}u).
\end{equation} In this case $\lambda=\frac{\widetilde{A}}{\sqrt{k_0 c_0}}$ and $D=0$.

The equation \eqref{fi1} for $\phi_1$ becomes 
\begin{equation}
    \phi_1(\xi)
    =
    \widetilde{E}
    + \tfrac{\widetilde{D}}{k_0}
      \int
      \exp{(\widetilde{A}\phi_1(\xi)})
      \exp\!\left(
        -\tfrac{k_0 c_0}{2\widetilde{A}^{2}}\int 
         \tfrac{\xi}{ \left(\int
      \exp{(-\widetilde{A}\phi_1(\xi)}) \,d\xi\right)^2 }\,\right)\,d\xi .
\end{equation}

 For $\phi_3$ we have \eqref{fi3} is equivalent to 
 \begin{equation}
       \phi_3(x) = u_1 k_0 \int \exp(\widetilde{A}\phi_3(x))\,dx,
 \end{equation} then
 \begin{equation}
      \phi^{'}_3(x)\exp(-\widetilde{A}\phi_3(x)) = u_1 k_0 ,
 \end{equation}
whose solution is given implicitly by
\begin{equation}
\phi_3(x)=-\frac{1}{\widetilde{A}}\ln(-u_1k_0\widetilde{A}x+\Hat{u}_1),
\end{equation} with $u_1$, $\hat{u}_1$ constants.

From \eqref{fi4} we have 

\begin{equation}\label{fi41}
   k_0 \int  \exp(-\widetilde{A}u)\,du= \frac{-2x}{\sqrt{Q}} ,
\end{equation}
then we obtain that solution $u$ is given explicitly by
\begin{equation}
u=-\frac{1}{\widetilde{A}}\ln\left(\frac{2\widetilde{A}x}{k_0 \sqrt{Q}}+\Hat{u}_2\right),
\end{equation} with $u_1$, $\hat{u}_2$ constants.
}

\subsection{{\color{black}{\color{black}Power} type $C(u)$ and $K(u)$}}

Finally, we consider PDE \eqref{eq} in the case where the coefficients \(C\) and \(K\) are given by
\begin{equation}
C(u)=\rho\,c_0\,(1+\beta u^p), \qquad
K(u)=k_0\,(1+\beta u^p),
\end{equation}
where $p$, \(\beta\), \(\rho\), \(c_0\), and \(k_0\) are nonnegative constants.
Coefficients of this type arise in free boundary problems that have been studied in several works, for instance, \cite{KuSiRa2020A,BoNaSeTa2022-Coam,BoCaNaSeTa2026}. 
Taking into account that
\[
\frac{C(u)}{K(u)}=\frac{\rho c_0}{k_0}=\alpha,
\]
and in accordance with Case~2 described in Sections~2 and~3, equation \eqref{eq}
admits the symmetry generators
\(\overline{X}_1\), \(\overline{X}_2\), \(\overline{X}_3\),
\(\overline{X}_4\), and \(\overline{X}_5\).
The corresponding invariant solutions, summarized in
Table~\ref{tab:invariant_solutions_case2}, are given by
\begin{itemize}
    \item $u$ is  implicitly defined by \eqref{psi1_eq}, which is equivalent to
    \begin{equation}
        u+\beta \frac{u^{p+1}}{p+1}=\frac{1}{\sqrt{t}}\left( \frac{ax}{t}+b\right) \exp\left(-\tfrac{\alpha x^2}{4t} \right), \quad a,b \in \mathbb{R}.
    \end{equation}
    \item From \eqref{psi2} we have that $u=\psi_2(\xi)$, $\xi=x/\sqrt{t}$ must be a solution of the integral equation given by
    \begin{equation}
        \psi_2(\xi)= \widetilde{E} +\widetilde{D}\int \frac{1}{1+\beta\psi_2^p(\xi)}\exp\left(-\frac{\alpha\xi^2}{4}\right) d\xi.
    \end{equation}
    \item $u$ is implicitly defined by \eqref{u-X1raya}, which is equivalent to
    \begin{equation}
        u+\beta\frac{u^{p+1}}{p+1}=\frac{a}{\sqrt{t}}\exp\left(-\frac{\alpha x^2}{4t}\right) +b, \quad a,b \in \mathbb{R}.
    \end{equation}
    \item $u \equiv u_0$
    \item From \eqref{psi5} we have $u=\psi_5(x)$ must be solution of
    \begin{equation}
        \psi_5(x)=a\int \frac{dx}{1+\beta\psi_5^p(x)}+b, \quad a,b \in \mathbb{R}.
    \end{equation}
\end{itemize} 
{\color{black} For the particular case $p=1$, that is, when $C$ and $K$ are linear functions of $u$, we obtain the invariant solutions associated with the generators 
$\overline{X}_1$, $\overline{X}_2$, $\overline{X}_3$, $\overline{X}_4$, and $\overline{X}_5$. }

\begin{itemize}

\item The function $u$ is given by~\eqref{psi1_eq}, namely,
\begin{equation}
\left(u + \frac{1}{\beta}\right)^2
=
\frac{1}{\sqrt{t}}
\left( \frac{\widetilde{a} x}{t} + \widetilde{b} \right)
\exp\left(-\frac{\alpha x^2}{4t}\right)+{\color{black}\frac{1}{\beta^{2}}},
\qquad \widetilde{a},\widetilde{b} \in \mathbb{R}.
\end{equation}

\item In this case, $u = \psi_2(\xi)$, where the similarity variable is 
$\xi = \dfrac{x}{\sqrt{t}}$. The function $\psi_2$ satisfies the integral equation~\eqref{psi2}, which is equivalent to
\begin{equation}
(1 + \beta \psi_2(\xi))^2
=
\hat{D}\,\erf\left(\frac{\sqrt{\alpha}\,\xi}{2} \right)+ \hat{E},
\qquad \hat{D}, \hat{E} \in \mathbb{R},
\end{equation}
where $\erf$ denotes the error function defined by
\begin{equation}
\erf(z) = \frac{2}{\sqrt{\pi}} \int_0^z e^{-s^2}\, ds.
\end{equation}
 \item $u \equiv u_0$
    \item From \eqref{psi5} we have $u=\psi_5(x)$ must satisfy
     \begin{equation}
(1+\beta \psi_5(x))^2
=
\hat{a}\,x + \hat{b},
\qquad \hat{a}, \hat{b} \in \mathbb{R}.
\end{equation}
\end{itemize}

\section{Comparison with the standard nonlinear diffusion equation}

In this section, we map the Lie point symmetries derived in section 2 for   equation \eqref{ec}
\begin{equation*}
C(u)u_t = (K(u)u_x)_x,
\end{equation*}
to those reported in \cite{BlKu1989, Ov1959}  for the standard nonlinear diffusion equation
\begin{equation} \label{eqD}
w_t = (\mathcal{D}(w)w_x)_x,
\end{equation}
where $\mathcal{D}(w)$ is a variable diffusion coefficient.

It is well known that, that under a Kirchhoff-type transformation
\begin{equation} \label{transf}
w=\mathcal{W}(u)=\int_{u_r}^{u} C(\tilde{u})\, d\tilde{u},
\end{equation}
where $u_r$ denotes a reference temperature, equation \eqref{ec} is transformed into equation \eqref{eqD}, with diffusion coefficient
\begin{equation} \label{D}
\widetilde{D}(w)=\frac{K(u)}{C(u)}.
\end{equation}

For arbitrary $\mathcal{D}(w)$, the infinitesimal generators of equation \eqref{eqD} are:
\begin{equation*}
\begin{aligned}
    X_1^w &= \frac{x}{2}\pder{x} + t\pder{t}, \\
    X_2^w &= \pder{x}, \\
    X_3^w &= \pder{t}.
\end{aligned}
\end{equation*}
If $\mathcal{D}(w) = \lambda (w+\kappa)^{\nu}$, an  infinitesimal generator is added, given by 
\begin{equation*}
    X_4^w = x\pder{x} + \frac{2}{\nu}(w+\kappa)\pder{w}.
\end{equation*}
In the case where $\mathcal{D}(w) = \lambda e^{\nu w}$, this becomes to 
\begin{equation*}
    X_4^w = x\pder{x} + \frac{2}{\nu}\pder{w}.
\end{equation*}
When $\mathcal{D}(w) = \lambda(w + \kappa)^{-\frac{4}{3}}$, a fifth generator is added, namely
\begin{equation*}
    X_5^w = x^2 \pder{x} - 3x(w+\kappa)\pder{w}.
\end{equation*}
The relations \eqref{transf} and \eqref{D} yield a one-to-one correspondence between the generators $X_i$ obtained in Section 2 and $X_i^w$, and between the corresponding nonlinear coefficients $C(u), K(u)$ and $\mathcal{D}(w)$, as summarized in the following scheme:
\begin{itemize}
    \item[(i)] For arbitrary functions $K(u)$ and $C(u)$ such that $\frac{K(u)}{C(u)}$ is non-constant (i.e., $\mathcal{D}(w)$ is an arbitrary non-constant coefficient), the following one-to-one correspondences hold:
    \begin{align*}
        X_1^w &\longleftrightarrow X_1, \\
        X_2^w &\longleftrightarrow X_2, \\
        X_3^w &\longleftrightarrow X_3.
    \end{align*}

    \item[(ii)] When $C(u)$ and $K(u)$ satisfy the condition
    \[
    C(u) = E K(u) \left(B \int K(u) \, du+D\right)^{\frac{1}{B}},
    \]
    the following subcases arise:
    \begin{itemize}
        \item[(a)] If $B \neq -1$ and $B \neq 0$, then $\mathcal{D}(w)$ adopts the power-law form
        \[
        \mathcal{D}(w) = \lambda (w + \kappa)^{\nu},
        \]
   with      $\nu = -\frac{1}{B+1}$, and $\lambda = \left(E^B (B+1)\right)^{-\frac{1}{B+1}}$. In this case, the additional generator $X_4$ corresponds to $X_4^w$ through the following chain of equalities:
        \[
        X_4 = x\pder{x} - 2\frac{B \int K(u) \, du+D}{K(u)}\pder{u}
            = x\pder{x} + \frac{2}{\nu}(w + \kappa)\pder{w}
            = X_4^w.
        \]
    \item[(b)] If $B=-1$, then $\mathcal{D}(w)$ takes the exponential form
    \[
    \mathcal{D}(w) = \lambda e^{\nu w},
    \]
    with $\nu = -\frac{1}{E}$ and $\lambda$ constant. The corresponding generator is given by
    \[
    X_4 = x\pder{x} + 2\frac{ \int K(u) \, du-D}{K(u)}\pder{u}=x\pder{x} + \frac{2}{\nu}\pder{w} = X_4^w.
    \]
\end{itemize}
\item[(iii)] If $B=0$, the relation between  $C(u)$ and $K(u)$ reduces to
\[
C(u) = E K(u) \exp\left(\frac{1}{D} \int K(u) \,du\right),
\]
which yields $\mathcal{D}(w) = D(w+\kappa)^{-1}$, corresponding to $\nu = -1$. In this case, we have
\[
X_4 = x\pder{x} - 2\frac{D}{K(u)}\pder{u}
    = x\pder{x} - 2(w+\kappa)\pder{w}
    = X_4^w.
\]
\item[(iv)] If
\[
C(u) = \frac{NK(u)}{\left(4M-\int K(u)\,du\right)^4},
\]
then $\mathcal{D}(w)$ adopts the form
\[
\mathcal{D}(w) = \lambda(w + \kappa)^{-\frac{4}{3}},
\]
with $\lambda = 3^{-\frac{4}{3}} N^{\frac{1}{3}}$, and the fifth generator reads
\[
X_5 = x^2\pder{x} - x\left(\frac{4M-\int K(u)\,du}{K(u)}\right)\pder{u}
    = x^2\pder{x} - 3x(w + \kappa)\pder{w}
    = X_5^w.
\]
\end{itemize}



Although equation \eqref{ec} can be reduced to equation \eqref{eqD}, whose Lie symmetry analysis is already available in the literature, in this paper we choose to work with the original formulation \eqref{ec}, since it preserves physical information that may be lost when the model is rewritten in terms of a single coefficient $\mathcal{D}$.

In particular, preserving the explicit decoupling of the specific heat $C$ and the thermal conductivity $K$ is essential for tracking the temporal evolution of the temperature field $u(x,t)$, thereby maintaining the physical traceability of the underlying heat storage and heat conduction mechanisms.

In materials with strongly temperature-dependent thermal properties, the specific heat $C(u)$ and the thermal conductivity $K(u)$ may vary at different, and sometimes even opposite, rates as the temperature changes. In such situations, the transient thermal behavior is determined by the local balance between thermal storage and heat transport mechanisms. The evolution of the temperature profile therefore depends on the relative influence of $C(u)$ and $K(u)$ at each point in space and time. By retaining the original formulation, this interplay remains explicitly represented, providing physical insight into the governing thermal processes that would otherwise be hidden when the model is recast in terms of a single effective diffusion coefficient $\mathcal{D}(w)$.


A further difficulty arises when boundary conditions are incorporated. In this case, transforming a problem governed by \eqref{ec} into an equivalent problem governed by \eqref{eqD} may introduce highly nonlinear implicit boundary conditions. To illustrate this issue, consider the following moving boundary problem:
\[
\begin{cases}
C(u),u_t=(K(u),u_x)_x, & 0<x<s(t), \quad t>0, \\
K(u(0,t)),u_x(0,t)=h_0(t)\bigl(u(0,t)-f_0(t)\bigr), & t>0,\\
u(s(t),t)=g_0(t), & t>0.
\end{cases}
\]
where $s(t)$ is a prescribed boundary. Applying the transformation \eqref{transf}, this problem becomes equivalent to
\[
\begin{cases}
w_t=(\mathcal{D}(w),w_x)_x, & 0<x<s(t), \quad t>0,\\
\mathcal{D}(w(0,t)),w_x(0,t)=h_0(t)\bigl(\mathcal{W}^{-1}(w(0,t))-f_0(t)\bigr), & t>0,\\
w(s(t),t)=\mathcal{W}(g_0(t)), & t>0.
\end{cases}
\]

It is readily observed that the boundary conditions in the transformed formulation become considerably more involved. In particular, the Robin condition is no longer expressed directly in terms of the physical temperature, but instead contains the inverse transformation $\mathcal{W}^{-1}$. Consequently, even when the original boundary conditions are simple and physically transparent, their counterparts in the diffusion formulation may acquire a highly nonlinear implicit structure. This example illustrates that preserving the original formulation not only retains the distinct physical roles of $C$ and $K$, but may also lead to a more natural and tractable description of the boundary conditions.

\section{Conclusion}
{\color{black}In this work, we have investigated a class of nonlinear heat-diffusion equations with temperature-dependent coefficients, focusing on the conditions under which the equations admit nontrivial Lie symmetry generators. By systematically analyzing the forms  $C(u)$ and $K(u)$, we identified specific functional dependencies that allow the existence of one-parameter Lie groups of transformations. Using these symmetries, invariant solutions were constructed for three representative cases found in the literature: power-type $C(u)$ with constant $K(u)$, coefficients satisfying the Storm's condition, and {\color{black}power}-type coefficients.

The results provide a clear framework to determine which forms of $C(u)$ and $K(u)$ yield solvable or partially solvable models via symmetry methods. These invariant solutions can serve as benchmarks for numerical simulations and as analytical tools for studying phase-change problems with variable thermal properties. Future work may extend this approach to  equations including source terms {\color{black} or free boundary problems}}

\section*{Acknowledgements}
{\color{black}This work was supported by Projects O06-26CI2100 and O06-26CI2101 Universidad Austral, Rosario, Argentina, Proyect PIP-CONICET 11220220100532CO, Argentina.
}
\bibliographystyle{plain} 
\bibliography{Biblio}

\end{document}